\documentclass[11pt]{amsart}
\usepackage{amsmath,amsthm,amsfonts,amssymb,latexsym,epsfig}

\newtheorem{theorem}{Theorem}[section]

\newtheorem{cor}{Corollary}[section]

\numberwithin{equation}{section}

\theoremstyle{definition}

\theoremstyle{remark}

\begin{document}
\title{On $l^p$ norms of factorable matrices}
\author{Peng Gao}
\address{Department of Mathematics, School of Mathematics and System Sciences, Beijing University of Aeronautics and Astronautics, P. R. China}
\email{penggao@buaa.edu.cn}


\begin{abstract}
  We study $l^p$ operator norms of factorable matrices and related results. We give applications to $l^p$ operator norms of weighted mean matrices and Copson's inequalities. We also apply the method in this paper to study the best constant in an inequality of Hardy, Littlewood and  P\'{o}lya.
\end{abstract}

\maketitle
\section{Introduction}
\label{sec 1} \setcounter{equation}{0}

  Suppose throughout the paper that $p\neq 0, \frac{1}{p}+\frac{1}{q}=1$.
  When $p \geq 1$, let $l^p$ be the Banach space of all complex sequences ${\bf x}=(x_n)_{n \geq 1}$ with norm
\begin{equation*}
   ||{\bf x}||_p: =(\sum_{n=1}^{\infty}|x_n|^p)^{1/p} < \infty.
\end{equation*}

    Let $C=(c_{n,k})$ be a matrix acting on the $l^p$ space. The $l^{p}$ operator norm of $C$ is
   defined as
\begin{equation*}
\label{02}
    ||C||_{p,p}=\sup_{||{\bf x}||_p = 1}\Big | \Big |C \cdot {\bf x}\Big | \Big |_p.
\end{equation*}

   It follows that for any $U_p \geq ||C||_{p,p}$ and any ${\bf x} \in
   l^p$,
\begin{equation*}
\label{01}
    \sum^{\infty}_{n=1} \Big{|}\sum^{\infty}_{k=1}c_{n,k}x_k
    \Big{|}^p \leq U_p \sum^{\infty}_{n=1}|x_n|^p.
\end{equation*}

   A prototype of the above inequality is the celebrated
   Hardy inequality (\cite[Theorem 326]{HLP}), which asserts that for $p>1$, the Ces\'aro matrix
    operator $C$, given by $c_{n,k}=1/n , 1 \leq k \leq n$ and $0$
    otherwise, has norm $p/(p-1)$.

    We say a matrix $A=(a_{n,k})$ is a lower triangular matrix if $a_{n,k}=0$ for $n<k$ and $A$ is a factorable matrix if it is a lower triangular matrix and in addition if $a_{n,k} = a_nb_k$ when $1 \leq k \leq n$. We say a factorable matrix $A$ is a weighted mean matrix if its entries satisfy:
\begin{equation}
\label{021}
    a_{n,k}=\lambda_k/\Lambda_n,  \quad 1 \leq k \leq
    n,
\end{equation}
   where in this paper, we let $(\lambda_n)_{n \geq 1}$ be a positive sequence and we assume all infinite sums converge. For two integers $N \geq n \geq 1$, we define
\begin{align}
\label{1.00}
  \Lambda_n=\sum^n_{k=1}\lambda_k,  \quad  \Lambda^*_{n}=\sum^{\infty}_{k=n}\lambda_k.
\end{align}
   Without loss of generality, we assume that the sequences ${\bf x}$ are positive in this paper and we further denote
\begin{align}
\label{1.01}
 & A_n=\sum^{n}_{k=1}\frac {\lambda_kx_k}{\Lambda_n},  \quad A^T_{n}=\lambda_n\sum^{\infty}_{k=n}\frac {x_k}{\Lambda_k},  \quad  A^*_{n}=\sum^{\infty}_{k=n}\frac {\lambda_kx_k}{\Lambda^*_{n}},   \quad {A^*_{n}}^T=\lambda_n\sum^{n}_{k=1}\frac {x_k}{\Lambda^*_k}.  
\end{align}

    Hardy's inequality motivates one to
    study the $l^{p}$ operator norm of weighted mean matrices. Recently, the author \cite{G7} proved the following result:
\begin{theorem}\cite[Theorem 3.1]{G7}
\label{thm3.1}
  Let $p>1$ be fixed. Let $A$ be a weighted mean matrix given by
    \eqref{021}. We have
    $||A||_{p,p} \leq p/(p-L)$ if for any integer $n \geq 1$, there exists a positive constant
    $0<L<p$ such that
\begin{align*}
    \sum^n_{k=1}\frac {\lambda_k}{\Lambda_n}\prod^{n}_{i=k}\Big (\frac {\Lambda_{
    i+1}/\lambda_{i+1}-L/p}{\Lambda_i/\lambda_i} \Big
    )^{1/(p-1)} \leq
     \frac {p}{p-L}.
\end{align*}
\end{theorem}

    The above theorem implies the following result in \cite{G5}:
\begin{cor}\cite[Theorem 1.2]{G5}
\label{thm03}
    Let $p>1$ be fixed. Let $A$ be a weighted mean matrix given by
    \eqref{021}. We have
    $||A||_{p,p} \leq p/(p-L)$ if for any integer $n \geq 1$, there exists a positive constant
    $0<L<p$ such that
\begin{equation*}
\label{024}
    \frac {\Lambda_{n+1}}{\lambda_{n+1}} \leq \frac
    {\Lambda_n}{\lambda_n}  \Big (1- \frac
    {L\lambda_n}{p\Lambda_n} \Big )^{1-p}+\frac {L}{p}.
\end{equation*}
\end{cor}

   The above corollary further implies the following well-known result of Cartlidge \cite{Car}:
\begin{cor}
\label{thm02}
    Let $p>1$ be fixed. Let $A$ be a weighted mean matrix given by
    \eqref{021}. We have
    $||A||_{p,p} \leq p/(p-L)$ if
\begin{align}
\label{022}
    L=\sup_n\Big(\frac {\Lambda_{n+1}}{\lambda_{n+1}}-\frac
    {\Lambda_n}{\lambda_n}\Big) < p.
\end{align}
\end{cor}

   We note that by slightly modifying the proof of Theorem \ref{thm3.1} in \cite{G7}, one obtains the following analogue of Theorem \ref{thm3.1} concerning the $l^p$ norms of factorable matrices:
\begin{theorem}
\label{thm3.1'}
  Let $p>1$ be fixed. Let $A=(a_{n,k})$ be a factorable matrix satisfying $a_{n,k}=a_nb_k$ when $1 \leq k \leq n$ and $a_n>0, b_n>0, a_1=b_1$. Let $0<L<p$ and $\lambda_p=(1-L/p)^p$. Then
    $||A||_{p,p} \leq p/(p-L)$ if for an integer $n \geq 1$, the following condition is satisfied:
\begin{align*}
  \sum^n_{k=1}\frac {b_k}{a_n}\prod^{n}_{i=k}\Big (\frac {a_{
    i}/b_{i+1}+1-L/p}{a_i/b_i} \Big
    )^{1/(p-1)} \leq
     \frac {p}{p-L}.
\end{align*}
\end{theorem}
   There is a need to study the $l^p$ norms and related results of factorable matrices besides the class of weighted mean matrices, as they have many applications. For example, the following two inequalities are related to the $l^p$ norms of the corresponding factorable matrices:
\begin{align}
\label{3.4'}
  \sum^{\infty}_{n=1}\left ( \lambda^{1/p}_n\Lambda^{-c/p}_n \sum^{n}_{k=1}\lambda^{1-1/p}_k \Lambda^{-(1-c/p)}_k x_k \right
   )^p & \leq \left ( \frac {p}{c-1}\right )^p
   \sum^{\infty}_{n=1}x^p_n, \quad 1<c \leq p; \\
\label{2.00}
   \sum^{\infty}_{n=1}\left( \frac {\Lambda^{\alpha}_n-\Lambda^{\alpha}_{n-1}}{\lambda^{1-1/p}_n\Lambda^{\alpha}_n}\sum^{n}_{k=1}
   \lambda^{1-1/p}_kx_k  \right )^p
     & \leq \left( \frac {\alpha p}{p-1} \right
     )^p\sum^{\infty}_{n=1}x^p_n, \quad p>1, \alpha>0.
\end{align}

 Inequality \eqref{3.4'} is equivalent to the following classical Copson inequalities \cite[Theorem 1.1, 2.1]{C}:
\begin{align}
\label{1.1'}
  & \sum^{\infty}_{n=1}\lambda_n\Lambda^{p-c}_nA^p_n \leq \left ( \frac {p}{c-1}\right )^p
   \sum^{\infty}_{n=1}\lambda_n\Lambda^{p-c}_nx^p_n, \ 1<c \leq p
   ; \\
  \label{1.2'}
  & \sum^{\infty}_{n=1}\lambda_n\Lambda^{p-c}_n\left ( \frac 1{\Lambda_n}\sum^{\infty}_{k=n}\lambda_k x_k \right
   )^p \leq \left ( \frac {p}{1-c}\right )^p
   \sum^{\infty}_{n=1}\lambda_n\Lambda^{p-c}_nx^p_n, \ 0 \leq c < 1 <p.
\end{align}
   In fact, it is shown in  \cite{G8} that 
   the above inequalities are equivalent to each other and the following inequalities of Leindler \cite[(1)]{L}:
\begin{align}
\label{1.3'}
 & \sum^{\infty}_{n=1}\lambda_n{\Lambda^*_n}^{p-c}\left (
 \frac 1{\Lambda^*_n}\sum^{n}_{k=1}\lambda_k x_k \right
   )^p \leq \left ( \frac {p}{1-c}\right )^p
   \sum^{\infty}_{n=1}\lambda_n{\Lambda^*_n}^{p-c}x^p_n, \ 0 \leq c
   <1
   ; \\
\label{1.4'}
  & \sum^{\infty}_{n=1}\lambda_n{\Lambda^*_n}^{p-c}{A^*_n}^p \leq \left ( \frac {p}{c-1}\right )^p
   \sum^{\infty}_{n=1}\lambda_n{\Lambda^*_n}^{p-c}x^p_n, \ 1 < c \leq
   p.
\end{align}
    The constants in \eqref{1.1'}-\eqref{1.4'} are best possible. 
  
    Bennett  
    \cite[p. 411]{B1} observed that
   inequality \eqref{1.1'} continues to hold for $c>p$ with
   constant $(p/(p-1))^p$. It is then natural to ask 
   whether inequality \eqref{1.1'} itself continues to hold for $c>p$. Note
   that in this case the constant $(p/(c-1))^p$ is best possible (see \cite{G8}). In \cite{G8}, the following extension of \eqref{1.1'} is given:
\begin{theorem}\cite[Theorem 2.1]{G8}
\label{thm1} Let $p>1$ be fixed. Let $c_p<0$ denote the unique
number satisfying
\begin{align*}
    \left (1+\frac {1-c_p}{p} \right
   )^{1-p}-\frac {1-c_p}{p} =0.
\end{align*}
   Then inequality \eqref{1.1'} holds for all $1< c \leq p-(p-1)c_q$. 
\end{theorem}
 
   It is shown in \cite{G8} that inequality \eqref{2.00} implies the following result of Bennett and Grosse-Erdmann \cite[Theorem 8]{BGE1}, which asserts that for $p \geq 1, \alpha \geq 1$,
\begin{align}
\label{1.5}
  \sum^{\infty}_{n=1}\lambda_n\left (\sum^{\infty}_{k=n}\Lambda^{\alpha }_{k}x_k \right )^p \leq (\alpha
  p+1)^p\sum^{\infty}_{n=1}\lambda_n\Lambda^{\alpha p}_{n}\left
  (\sum^{\infty}_{k=n}x_k \right )^p.
\end{align}
  Here the constant is best possible. It is conjectured \cite[p.
   579]{BGE1} that inequality \eqref{1.5} (resp. its reverse) remains valid with the same best possible
   constant when $p \geq 1, 0 < \alpha < 1$ (resp. $ -1/p < \alpha < 0$).
   In \cite{G8}, the following partial resolution on the above mentioned conjecture is given:
\begin{theorem}
\label{thm2} Inequality \eqref{2.00} and hence \eqref{1.5} is valid for $p>1, \alpha \geq
1-\frac {1}{2p}$.
\end{theorem}

   Motivated by the above example, it is our goal in this paper to study the $l^p$ norms of factorable matrices and related results. We prove in the next section the following
\begin{theorem}
\label{mainthm}
  Let $p>1$ be fixed. Let $A=(a_{n,k})$ be a factorable matrix satisfying $a_{n,k}=a_nb_k$ when $1 \leq k \leq n$ and $a_n>0, b_n>0, a_1=b_1$. Let $0<L<p$ and $\lambda_p=(1-L/p)^p$. Then
    $||A||_{p,p} \leq p/(p-L)$ if for any integer $n \geq 1$, the following condition is satisfied:
\begin{align}
\label{2.8}
 & (\lambda^{1-1/p}_p\frac {a_{n+1}}{b_{n+1}}+1-\lambda^{1-1/p}_p)^{\frac 1{p-1}}\left ( (\lambda^{1-1/p}_p\frac
    {a_n}{b_n}+1-\lambda_p-\lambda^{1-1/p}_p)^{\frac 1{p-1}}+(\frac {a_{n}}{b_{n+1}})^{\frac p{p-1}} \right ) \\
\leq & \left ( \frac {a_{n+1}}{b_{n+1}}\right )^{\frac p{p-1}}(\lambda^{1-1/p}_p\frac
    {a_n}{b_n}+1-\lambda_p-\lambda^{1-1/p}_p)^{\frac 1{p-1}}.  \nonumber
\end{align} 
\end{theorem}

  For the weighted mean matrices case, we shall not compare Theorem \ref{mainthm} with Theorem \ref{thm3.1} (or rather, Corollary \ref{thm03} as it is more practical in applications) in general here.  We only point out that when $p=2$, Theorem \ref{mainthm} implies the following
\begin{cor}
\label{cor1.3}
  Let $A$ be a weighted mean matrix given by
    \eqref{021} and let $0<L<2$. Then
    $||A||_{2,2} \leq 2/(2-L)$ if for any integer $n \geq 1$, the following condition is satisfied:
\begin{align}
\label{1.04}
 \frac {\Lambda_{n+1}}{\lambda_{n+1}}-\frac
    {\Lambda_n}{\lambda_n} \leq L+\frac {L^2}{4}\cdot \frac {1+\frac L2}{1-\frac L2}\cdot \frac 1{\frac {\Lambda_{n+1}}{\lambda_{n+1}}+\frac{L}{2}}.
\end{align}
\end{cor}

  On the other hand, the condition given in Corollary \ref{thm03} when $p=2$ becomes
\begin{align}
\label{1.05}
\frac {\Lambda_{n+1}}{\lambda_{n+1}}-\frac
    {\Lambda_n}{\lambda_n} \leq L+\frac {L^2}{4}\cdot \frac {1}{\frac
    {\Lambda_n}{\lambda_n}-\frac{L}{2}}.
\end{align}
  It is then easy to see that inequality \eqref{1.04} and \eqref{1.05} are not comparable in general.


    The proof of Theorem \ref{thm1} given in \cite{G8} is via the study of an equivalent inequality of \eqref{1.1'} by the duality principle. Here and in what follows, we shall refer to the duality principle as the fact that the norms of a bounded linear operator and its adjoint operator coincide (see for example, \cite[Chap. 4]{Rudin}) and that the $l^p$ and $l^q$ spaces (with $p>1$) are dual spaces to each other (see for example, \cite[Chap. 19]{Royden}). In Section \ref{sec 4}, using different approaches to $l^p$ norms of factorable matrices, we give two other proofs of Theorem \ref{thm1}. We also give two other proofs of Theorem \ref{thm2} as well.

    Now, we consider some examples that are closely related to the study of $l^p$ norms of factorable matrices. They can be regarded as another motivation for the paper.  
    
   We note first that the above mentioned result of Cartlidge has the following strengthened form (see \cite{G1}):
\begin{theorem}
\label{Cartimprovement}
 Let $p>1$ be fixed. Let $A$ be a weighted mean matrix. If \eqref{022}
is satisfied, then
\begin{align}
\label{1.40}
  \sum^{\infty}_{n=1}A^p_n \leq \left(\frac p{p-L} \right )\sum^{\infty}_{n=1}x_nA^{p-1}_n.
\end{align}
 In particular, $\| A \|_{p,p} \leq p/(p-L)$.
\end{theorem}
   As it is easy to show by H\"older's inequality (see \cite{G1}) that  \eqref{1.40} implies that
\begin{align}
\label{1.6}
  \sum^{\infty}_{n=1}A^p_n \leq \left(\frac p{p-L} \right )^p\sum^{\infty}_{n=1}x^p_n,
\end{align}
    the statement of Corollary \ref{thm02} therefore follows from that of Theorem \ref{Cartimprovement}. Similar to Theorem \ref{Cartimprovement}, we note the following
\begin{theorem}[{\cite[Theorem 6.3]{G8}}]
\label{thm3.20}
  Let $p>1$ be fixed. If
\begin{align}
\label{1.7}
L' = \sup_{n}\left (\frac {\Lambda^*_{n}}{\lambda_{n}}-\frac {\Lambda^*_{n+1}}{\lambda_{n+1}} \right )<p,
\end{align}
    then 
\begin{align}
\label{1.8'}
   \sum^{\infty}_{n=1}{A^*_n}^p \leq \left (\frac {p}{p-L} \right )
    \sum^{\infty}_{n=1}x_n {A^*_n}^{p-1}.
\end{align}
\end{theorem}

   One can show by the method of sinister transpose \cite[p. 408]{B1} that inequalities \eqref{1.40} and \eqref{1.8'} are equivalent.  It follows from the duality principle that inequality \eqref{1.6} is equivalent to the following (with $L<p/(p-1)$ )
\begin{align}
\label{1.07'}
   \sum^{\infty}_{n=1}({A^T_n})^p \leq \left ( \frac {p}{p-(p-1)L} \right )^p \sum^{\infty}_{n=1}x^p_n.
\end{align}
   It is then natural to expect the following inequality to hold under condition \eqref{022}:
\begin{align}
\label{1.07}
    \sum^{\infty}_{n=1}({A^T_n})^p \leq \left ( \frac {p}{p-(p-1)L} \right ) \sum^{\infty}_{n=1}x_n({A^T_n})^{p-1}.
\end{align}
  Again, it is easy to show that inequality \eqref{1.07} implies \eqref{1.07'}. By the method of sinister transpose, the above inequality is equivalent to the following one under condition \eqref{1.7}:
\begin{align}
\label{1.9'}
    \sum^{\infty}_{n=1}({{A^*_n}^T})^p \leq \left ( \frac {p}{p-(p-1)L} \right ) \sum^{\infty}_{n=1}x_n({{A^*_n}^T})^{p-1}.
\end{align}
   
    In Section \ref{sec 5}, we show that the expectation above is true by proving the following
\begin{theorem}
\label{thm1.4}
  Let $p>1$ be fixed. Let $A$ be a weighted mean matrix. If \eqref{022}
is satisfied, then inequality \eqref{1.07} is valid. If  \eqref{1.7} is satisfied, then inequality \eqref{1.9'} is valid.
\end{theorem}
   
    Also, motivated by Theorems \ref{Cartimprovement} and \ref{thm3.20}, we show in Section \ref{sec 5} the following
\begin{theorem}
\label{thm1.5} Let $p>1$ and let $c$ be a constant. Under the notations given in \eqref{1.00} and \eqref{1.01}, the following inequalities hold: 
\begin{align}
\label{1.8}
 \sum^{\infty}_{n=1} \lambda_n\Lambda^{p-c}_nA^{p}_{n} 
& \leq \frac {p}{c-1}\sum^{\infty}_{n=1}\lambda_n\Lambda^{p-c}_nx_nA^{p-1}_n, \quad 1 < c \leq p; \\
\label{1.90}
\sum^{\infty}_{n=1} \lambda_n\Lambda^{p-c}_n\left ( \frac 1{\Lambda_n}\sum^{\infty}_{k=n}\lambda_k x_k \right
   )^p
& \leq \frac {p}{1-c}\sum^{\infty}_{n=1}\lambda_n\Lambda^{p-c}_nx_n\left ( \frac 1{\Lambda_n}\sum^{\infty}_{k=n}\lambda_k x_k \right
   )^{p-1}, \quad 0 \leq c <1; \\
\label{1.10}
\sum^{\infty}_{n=1} \lambda_n{\Lambda^{*}_n}^{p-c}\left ( \frac 1{\Lambda^*_n}\sum^{n}_{k=1}\lambda_k x_k \right
   )^p
& \leq \frac {p}{1-c}\sum^{\infty}_{n=1}\lambda_n{\Lambda^{*}_n}^{p-c}x_n\left ( \frac 1{\Lambda^*_n}\sum^{n}_{k=1}\lambda_k x_k \right
   )^{p-1}, \quad 0\leq c<1; \\
\label{1.11}
\sum^{\infty}_{n=1}\lambda_n{\Lambda^*_n}^{p-c}{A^*_n}^p &  \leq \frac {p}{c-1}
  \sum^{\infty}_{n=1}\lambda_n{\Lambda^*_n}^{p-c}x_n{A^*_n}^{p-1}, \ 1 < c \leq
   p.
\end{align}
\end{theorem} 
   We point out here the case $c=p$ of inequality \eqref{1.8} is shown by Copson in \cite{C} and the case $c=0$ of inequality \eqref{1.10} is shown by Leindler in \cite{L}. It it easy to see that the above inequalities imply inequalities \eqref{1.1'}-\eqref{1.4'}.  Moreover, inequality \eqref{1.8} is equivalent to inequality \eqref{1.11} and inequality \eqref{1.90} is equivalent to inequality \eqref{1.10} by the method of sinister transpose.

  We end our introduction by considering an inequality of Hardy, Littlewood and  P\'{o}lya \cite[Theorem 345]{HLP}, which  asserts that the following inequality holds for $0<p<1$ with $C_p=p^p$:
\begin{align}
\label{1}
  \sum^{\infty}_{n=1}\Big( \frac 1{n} \sum^{\infty}_{k=n}x_k \Big
  )^p \geq C_p \sum^{\infty}_{n=1}x^p_n.
\end{align}
   It is noted in \cite{HLP} that the constant $C_p=p^p$ may not be best possible
  and the best constant $C_p=(p/(1-p))^p$ when $0< p \leq \frac 13$ was indeed obtained by Levin and Ste\v ckin in \cite[Theorem 61]{L&S}. In \cite{G9}, the following extension of the above result of Levin and Ste\v ckin is given:
\begin{theorem}\cite[Theorem 1]{G9}
\label{thm6}
  Inequality \eqref{1}
  holds with the best possible constant $C_p=(p/(1-p))^p$ for any $1/3<p<1/2$ satisfying
\begin{equation}
\label{1.4}
   2^{p/(1-p)}\Big (\Big(\frac {1-p}{p} \Big
)^{1/(1-p)}-\frac {1-p}{p}\Big ) -(1+\frac {3-1/p}2)^{1/(1-p)}
\geq 0.
\end{equation}
   In particular, inequality \eqref{1}
  holds for $0 < p  \leq 0.346$.
\end{theorem}

   We note that it is shown in \cite{G9} that inequality \eqref{1} with $C_p=(p/(1-p))^p$ is equivalent to the following one (note that $x_n>0$ in this paper):
\begin{align}
\label{1.9}
  \sum^{\infty}_{n=1}\Big(\sum^{n}_{k=1}\frac {x_k}{k} \Big
  )^{q} \leq  \Big ( \frac {p}{1-p} \Big )^q \sum^{\infty}_{n=1}x^{q}_n.
\end{align}
   As the corresponding matrix of the above inequality is a factorable matrix, one can study the above inequality using methods similar to those used in the study of $l^p$ norms of factorable matrices.

   In Section \ref{sec 3}, we study inequality \eqref{1} with $C_p=(p/(1-p))^p$  and \eqref{1.9}. We first prove the following
\begin{theorem}
\label{thm3.2}
  Let $1/3 \leq p < 1$. Define a sequence $( \mu_n )_{n \geq 1}$ as
\begin{align*}
 \mu_1=\left ( \frac {1-p}{p} \right )^p, \quad \mu_{n+1}=(n+1)^p\left (n^{p/(p-1)}\mu^{1/(1-p)}_n-1 \right )^{1-p} +\left ( \frac {1-p}{p} \right )^p, \quad n \geq 1.
\end{align*}
 If there exists an integer $n_0 \geq 1$ such that 
\begin{align*}
  \mu_{n_0} \geq \left ( \frac p{1-p} \right )^{1-p}n_0+\frac 12 \cdot \left ( \frac 1p-1 \right )^p,
\end{align*} 
  then inequality \eqref{1}
  holds with the best possible constant $C_p=(p/(1-p))^p$.
\end{theorem}
   
  If we take $n_0=1$ in Theorem \ref{thm3.2}, then $\mu_1 \geq (p/(1-p))^{1-p}+(1/p-1)^p/2$ implies $p \leq 1/3$. If we take $n_0=2$, then $\mu_2 \geq 2(p/(1-p))^{1-p}+(1/p-1)^p/2$ is precisely inequality \eqref{1.4}. Thus, Theorem \ref{thm3.2} gives an improvement of Theorem \ref{thm6}. Calculations show that inequality \eqref{1} holds for $p=0.35$ by taking $n_0=4$ in Theorem \ref{thm3.2}. One then expects that
inequality \eqref{1} holds for all $0 < p \leq 0.35$ but we shall not worry about these numerical values here. 

   Theorem \ref{thm3.2} is proved by studying inequality \eqref{1} directly. We then study inequality \eqref{1.9} to prove the following
\begin{theorem}
\label{thm1.9}
  Let $1/3 \leq p < 1$. Define a sequence $( \mu_n )_{n \geq 1}$ as
\begin{align*}
 \mu_1=0, \quad \mu_{n+1}= \left( n^{-p}+\mu^{1-p}_n \right )^{1/(1-p)}-\left (\frac 1p-1 \right )^{p/(p-1)}, \quad n \geq 1.
\end{align*}
 If there exists an integer $n_0 \geq 1$ and a constant $c$ such that 
\begin{align}
\label{1.29}
  & \mu_{n_0} \geq (1/p-1)^{1/(p-1)}(n_0+c), \quad c>-\frac 1{2p}, \\
  & \left (\frac 1p-1 \right )\cdot \frac 1{n_0}+(1+c\cdot \frac 1{n_0})^{1-p} -\left (1+ \left (c+\frac 1p \right )\cdot \frac 1{n_0} \right )^{1-p} \geq 0, \nonumber
\end{align} 
  then inequality \eqref{1}
  holds with the best possible constant $C_p=(p/(1-p))^p$.
\end{theorem}

    One verifies that the conditions given in \eqref{1.29} are satisfied when $p=0.35, c=-1.33542621$. Hence Theorem \ref{thm1.9} also implies that inequality \eqref{1} holds with $C_p=(p/(1-p))^p$ when $p=0.35$. Again we shall leave further explorations on the numerical values of $p$ that make the validity of inequality \eqref{1} via Theorem \ref{thm1.9} to the interested reader.

\section{Proof of Theorem \ref{mainthm}}
\label{sec 2} \setcounter{equation}{0}
   Let $p>1$ and let $( a_n )_{n \geq 1}, ( b_n )_{n \geq 1}$ be two sequences of positive sequences satisfying $a_1=b_1$. We seek for conditions on $a_n, b_n$,
    such that the following inequality holds with a positive constant $U_p$:
\begin{align}
\label{2.100}
   \sum^{\infty}_{n=1}\Big{(}\sum^{n}_{k=1} \frac {b_kx_k}{a_n}
   \Big{)}^p \leq U_p \sum^{\infty}_{n=1}x_n^p.
\end{align}
  
   For this, first note that in order for inequality \eqref{2.100} to be valid, it suffices to establish the validity of it with the
infinite sums replaced by finite sums from $1$ to $N$, where $N \geq 1$ is an arbitrary integer. 
  We define
\begin{equation*}
    \sum^{n}_{k=1} \frac {b_kx_k}{a_n}=y_n.
\end{equation*}
   This allows us to recast inequality \eqref{2.100} as
\begin{equation}
\label{2.20}
   \lambda_p\sum^{N}_{n=1}y^p_n \leq  \sum^{N}_{n=1}
   \Big(\frac {a_n}{b_n}y_n-\frac {a_{n-1}}{b_{n}}y_{n-1}\Big
   )^p,
\end{equation}
  where we set $a_{0}=y_0=0, \lambda_p=U^{-1}_p$ and we require that $y_n \geq a_{n-1}y_{n-1}/a_n \geq 0$.

  For any integer $n \geq 1$ and fixed constants $\alpha>0$, $\beta>0$, $y_{n} \geq 0, \mu_n \geq 0$, we consider the following function for $0 \leq x \leq
  \alpha y_{n}/\beta$,
\begin{align}
\label{2.2'}
   f(x)=\Big(\alpha y_n-\beta x\Big
   )^p+\mu_nx^p.
\end{align}
   It is readily checked that $f'(x_0)=0$ implies that
\begin{align*}
    \alpha y_n-\beta x_0=\Big (\frac
{\mu_n}{\beta} \Big )^{\frac 1{p-1}}x_0.
\end{align*}
  Solving this for $x_0$, we obtain
\begin{equation*}
    x_0=\frac {\alpha \beta^{\frac 1{p-1}} y_n}{\mu^{\frac 1{p-1}}_n+\beta^{\frac p{p-1}}}.
\end{equation*}
  As $\mu_n \geq 0$, we have $0 \leq x_0 \leq \alpha y_{n}/\beta$. It is also easy to check that $f''(x_0) \geq
  0$, thus for all $0 \leq x \leq 
  \alpha y_{n}/\beta$, we have
\begin{align*}
   f(x) \geq  f(x_0)=\frac {\alpha^p\mu_ny^p_n}{\Big ( \mu^{\frac 1{p-1}}_n+\beta^{\frac p{p-1}}\Big )^{p-1}}.
\end{align*}
   We now set $x=y_{n-1}, \alpha=\frac {a_n}{b_n}, \beta=\frac {a_{n-1}}{b_n}$ to see that for all $n \geq 2$,
we have
\begin{align*}
  \Big(\frac {a_n}{b_n}y_n-\frac {a_{n-1}}{b_{n}}y_{n-1}\Big
   )^p \geq \frac {(a_n/b_n)^p \mu_n}{\Big(\mu^{1/(p-1)}_n+(a_{n-1}/b_n)^{p/(p-1)} \Big )^{p-1}}y^p_n-\mu_ny^p_{n-1}.
\end{align*}
   As $a_0=y_0=0$, we note that the above inequality continues to hold when $n=1$ for $\mu_1=1$. Summing up from $n=1$ to $N$ (with $\mu_{N+1} \geq 0$ as well), we obtain
\begin{align*}
   \sum^N_{n=1}\Big(\frac {a_n}{b_n}y_n-\frac {a_{n-1}}{b_{n}}y_{n-1}\Big
   )^p \geq \sum^{N}_{n=1}\left ( \frac {(a_n/b_n)^p \mu_n}{\Big(\mu^{1/(p-1)}_n+(a_{n-1}/b_n)^{p/(p-1)} \Big )^{p-1}}-\mu_{n+1}\right )y^p_n.
\end{align*}

   For $1 \leq n \leq N-1$, we define the values of $\mu_n$ inductively as follows (note that $\mu_1=1$):
\begin{align}
\label{recurrent}
 \mu_{n+1}=\frac {(a_n/b_n)^p \mu_n}{\Big(\mu^{1/(p-1)}_n+(a_{n-1}/b_n)^{p/(p-1)} \Big )^{p-1}}-\lambda_p.
\end{align}
  It is easy to see that the above values of $\mu_n$ lead to inequality \eqref{2.20} provided that the condition
  $\mu_n \geq 0$ is satisfied. 

  To prove Theorem \ref{mainthm}, we set $\lambda_p=\Big (1-\frac {L}{p} \Big )^p$ in \eqref{recurrent}. We may assume $N \geq 3$ and we proceed inductively to see what conditions will be imposed on $a_n, b_n$ so that we can have $\mu_n \geq a a_{n-1}/b_{n-1}+b$ for $n \geq 2$ with constants $a,b$ to be specified later so that $\mu_n \geq 0$ is
  satisfied. First note that we must have $a \geq 0$ and the case 
  $n=2$ implies that $a+b \leq 1-\lambda_p$.  Suppose that $\mu_n \geq
  a a_{n-1}/b_{n-1}+b$ is satisfied for some $2 \leq n \leq N-1$, then we have
\begin{align*}
   \mu_{n+1} \geq \frac {(a_n/b_n)^p (aa_{n-1}/b_{n-1}+b)}{\Big((a a_{n-1}/b_{n-1}+b)^{1/(p-1)}+(a_{n-1}/b_n)^{p/(p-1)} \Big )^{p-1}}-\lambda_p.
\end{align*} 
   Thus, it suffices to find conditions on $a_n, b_n$ such that
\begin{align}
\label{2.07}
 \frac {(a_n/b_n)^p (a a_{n-1}/b_{n-1}+b)}{\Big((a a_{n-1}/b_{n-1}+b)^{1/(p-1)}+(a_{n-1}/b_n)^{p/(p-1)} \Big )^{p-1}}-\lambda_p \geq a\frac {a_n}{b_n}+b.
\end{align}
   To proceed further, we restrict ourself to the weighted mean matrices case by setting $a_n=\Lambda_n, b_n=\lambda_n$ with $\Lambda_n=\sum^n_{i=1}\lambda_i$. We then deduce that $a_{n-1}/b_n=a_{n}/b_{n}-1$ and we recast the above inequality by letting $y=a_n/b_n, x=a_{n-1}/b_{n-1}$ as
\begin{align}
\label{2.5}
  (ay+b+\lambda_p)^{1/(p-1)}\left ( (ax+b)^{1/(p-1)}+(y-1)^{p/(p-1)} \right )-y^{p/(p-1)}(ax+b)^{1/(p-1)} \leq 0.
\end{align} 
   To find the optimal values of $a, b$, we assume that $y-x=L+O(1/y)$ and consider the case $y \rightarrow +\infty$. We use Taylor expansion to determine the optimal values of $a, b$. As it is easy to show that the $O(1/y)$ term plays no role in this process, we may simply set $x=y-L$ in \eqref{2.5} to cast the left-hand side expression of \eqref{2.5} as
\begin{align}
\label{2.6}
 &  (ay)^{1/(p-1)}\left (1+\frac{b+\lambda_p}{ay} \right )^{1/(p-1)}\left (  (ay)^{1/(p-1)}(1+\frac {b-aL}{ay})^{1/(p-1)}+y^{p/(p-1)}(1-\frac 1y)^{p/(p-1)} \right ) \\
&-y^{p/(p-1)}(ay)^{1/(p-1)}(1+\frac {b-aL}{ay})^{1/(p-1)}. \nonumber
\end{align} 
   We consider the Taylor expansions of the following terms of the above expression 
\begin{align}
\label{2.7}
 \left (1+\frac{b+\lambda_p}{ay} \right )^{1/(p-1)}, \quad (1+\frac {b-aL}{ay})^{1/(p-1)}, \quad (1-\frac 1y)^{p/(p-1)}
\end{align} 
   to see that the leading term in \eqref{2.6} is
\begin{align*}
  \left (\lambda_p+a(L-p)+(p-1)a^{p/(p-1)}\right )\frac {a^{(2-p)/(p-1)}y^{2/(p-1)}}{p-1}.
\end{align*}
   It is easy to see that $\lambda_p+a(L-p)+(p-1)a^{p/(p-1)} \geq 0$ for $a \geq 0$ with equality holding if and only if $a=\lambda^{1-1/p}_p$. This implies that we must take $a=\lambda^{1-1/p}_p$. To determine the value of $b$, we proceed as above by considering the Taylor expansions in \eqref{2.7} to see that the coefficients involving with $b$ of the second leading term $y^{2/(p-1)-1}$ of \eqref{2.6} is
\begin{align*}
 \frac {a^{1/(p-1)-2}}{p-1}\left (2a^{p/(p-1)}-\frac {pa}{p-1}+\frac {(2-p)(\lambda_p+aL)}{p-1} \right )b=- \frac {a^{1/(p-1)-1}L}{p-1}b.
\end{align*}
  This implies that we should take the value of $b$ to be as large as possible. As $a+b \leq 1-\lambda_p$, we see that we need to take $b=1-\lambda_p-a=1-\lambda_p-\lambda^{1-1/p}_p$. 

Now, by setting $a=\lambda^{1-1/p}_p, b=1-\lambda_p-\lambda^{1-1/p}_p$ in \eqref{2.07} and by letting $n \rightarrow n+1$ there, we see that inequality \eqref{2.07} coincides with inequality \eqref{2.8} and this completes the proof of Theorem \ref{mainthm}.
\section{Copson's inequalities and a related result}
\label{sec 4} \setcounter{equation}{0}

   We first give two proofs of Theorem \ref{thm1}. Instead of \eqref{1.1'}, we study inequality \eqref{3.4'} here. For the first proof, we set $a_n=\lambda^{-1/p}_n\Lambda^{c/p}_n, b_n=\lambda^{1-1/p}_n \Lambda^{-(1-c/p)}_n, U_p=(p/(c-1))^p$ in \eqref{2.100} to see that, from our discussion in the previous section, in order for inequality \eqref{3.4'} to hold, it suffices to have a non-negative sequence $( \mu_n )_{n \geq 1}$ defined inductively as in \eqref{recurrent} with $\mu_1=1$. Explicitly, if we set $y_n=\lambda_n/\Lambda_n$ for $n \geq 1$ and $y_0=1$, the recurrence relation \eqref{recurrent} becomes
\begin{align}
\label{3.4}
 \mu_{n+1}=\frac {y^{-p}_n \mu_n}{\Big(\mu^{1/(p-1)}_n+y^{-1/(p-1)}_{n-1}y^{-1}_n(1-y_n)^{(c-1)/(p-1)} \Big )^{p-1}}-\left(\frac {c-1}p \right )^p.
\end{align}

   If there exists values of $c > p$ such that $\mu_n \geq (a y_{n-1})^{-1}$ for all $n \geq 2$, where $a$ is a positive constant to be determined in what follows, then these values of $c$ make inequality \eqref{1.1'} valid. To simplify the notations, we set $y=y_n$ and substitute the lower bound of $\mu_n$ in \eqref{3.4} to see that it suffices to find values of $c > p$ such that for $0 \leq y \leq 1$ (note that $y=1$ corresponds to the case $n=2$),
\begin{align}
\label{3.5}
   1+a\left(\frac {c-1}p \right )^p y \leq (a^{-1/(p-1)}y+(1-y)^{(c-1)/(p-1)} )^{1-p}.
\end{align}
   By considering the coefficient of $y$ of the Taylor expansions of the expressions on both sides of the above inequality, we see that we must have
\begin{align*}
  a\left(\frac {c-1}p \right )^p+(p-1)a^{-1/(p-1)}-(c-1) \leq 0.
\end{align*}
   The left-hand side expression above, when regarded as a function of $a$, is minimized at $a=((c-1)/p)^{1-p}$ with value $0$. Thus, we need to take $a=((c-1)/p)^{1-p}$. We can thus simplify inequality \eqref{3.5} as
\begin{align}
\label{3.6'}
  1+\left(\frac {c-1}p \right ) y \leq (\frac {c-1}p y+(1-y)^{(c-1)/(p-1)} )^{1-p}.
\end{align}
   We point out that we can identify the above inequality with inequality (2.1) in \cite{G8} by setting $c=(p-c)/(p-1), p=q, x=y$ there. We then deduce easily the statement of Theorem \ref{thm1} by Lemma 2.1 of \cite{G8}.

   We now give another proof of Theorem \ref{thm1}. Note that the discussion on the $l^p$ norms of weighted mean matrices via the duality principle in Section 5 of \cite{G5} carries over to factorable matrices as well, once one replaces $\Lambda_n, \lambda_n$ by $a_n, b_n$. In particular, the last equation on \cite[p. 843]{G5} implies that in order for inequality \eqref{2.100} to hold, it suffices to define a positive sequence $( \mu_n )_{n \geq 1}$:
\begin{align}
\label{3.7}
 \mu_1=\mu_{n+1}-\frac {1}{\Big
(\Big ( a_n/b_{n} \Big )^{q/(q-1)}\mu_n^{-1/(q-1)}-1
\Big )^{q-1}}\Big ( \frac {a_n}{b_{n+1}}\Big )^q =U^{-q/p}_p
\end{align}
  satisfying $\mu_n < (a_n/b_n)^q$.

  We now set $a_n=\lambda^{-1/p}_n\Lambda^{c/p}_n, b_n=\lambda^{1-1/p}_n \Lambda^{-(1-c/p)}_n, U_p=(p/(c-1))^p$ in \eqref{3.7} to see that inequality \eqref{3.4'} and hence inequality \eqref{1.1'} follows provided that we can define a positive sequence $( \mu_n )_{n \geq 1}$:
\begin{align*}
 \mu_1=\mu_{n+1}-\frac {(\frac {\lambda_n}{\Lambda_n})(\frac {\Lambda_n}{\Lambda_{n+1}})^{(c-1)/(p-1)}(\frac {\Lambda_{n+1}}{\lambda_{n+1}})}{\Big
(\mu_n^{-1/(q-1)}-\Big ( \lambda_n/\Lambda_n \Big )^{q/(q-1)}
\Big )^{q-1}}=\left(\frac {c-1}p \right )^{q}
\end{align*}
  satisfying $\mu_n < (\Lambda_n/\lambda_n)^q$.

  If there exists values of $c > p$ such that for $n \geq 1$, 
\begin{align}
\label{3.9}
  \mu_n \leq \frac {\Lambda_n}{\lambda_n}\left (\frac {\lambda_n}{\Lambda_n}+a \right )^{1-q},
\end{align}
  where $a$ is a positive constant to be determined in what follows, then these values of $c$ make inequality \eqref{1.1'} valid. To simplify the notations, we set $y=\lambda_{n+1}/\Lambda_{n+1}$ and substitute the upper bound of $\mu_n$ in \eqref{3.9} to see that it suffices to find values of $c > p$ such that for $0 \leq y \leq 1$ (note that $y=1$ corresponds to the case $n=1$),
\begin{align*}
  (1-y)^{\frac {c-1}{p-1}}+a^{q-1}\left ( \frac {c-1}{p} \right )^qy \leq (1+\frac {y}{a})^{1-q}.
\end{align*}
   Once again by considering the Taylor expansions, we see that we must take $a=p/(c-1)$ and the above inequality then becomes inequality \eqref{3.6'} and this gives another proof of Theorem \ref{thm1}.


     We now give two proofs of Theorem \ref{thm2}. We study inequality \eqref{2.00} here. For the first proof, we set 
\begin{align}
\label{3.13'}
 a_n=\frac {\lambda^{1-1/p}_n\Lambda^{\alpha}_n}{\Lambda^{\alpha}_n-\Lambda^{\alpha}_{n-1}}, \quad b_n=\lambda^{1-1/p}_n, \quad U_p= \left( \frac {\alpha p}{p-1} \right )^p
\end{align} 
   in \eqref{3.7} to see that inequality \eqref{2.00} is valid provided that we can define a positive sequence $( \mu_n )_{n \geq 1}$:
\begin{align*}
 \mu_1=\mu_{n+1}-\frac {(\frac {\lambda_n}{\Lambda_n})(\frac {\Lambda_{n}}{\lambda_{n+1}})}{\Big
(\mu_n^{-1/(q-1)}-\left ( 1- \left ( \Lambda_{n-1}/\Lambda_n \right )^{\alpha} \right )^{q/(q-1)}
\Big )^{q-1}}=\left(\frac {p-1}{\alpha p} \right )^{q}
\end{align*}
  satisfying $\mu_n < \left ( 1- \left ( \Lambda_{n-1}/\Lambda_n \right )^{\alpha} \right )^{-q}$.

  If there exists values of $0< \alpha < 1$ such that for $n \geq 1$,
\begin{align}
\label{3.14}
  \mu_n \leq \left ( \left ( 1- \left ( \Lambda_{n-1}/\Lambda_n \right )^{\alpha} \right )^{q/(q-1)}+\left (a\frac {\lambda_n}{\Lambda_n} \right )^{1/(q-1)} \right )^{1-q},
\end{align}
  where $a$ is a positive constant to be determined in what follows, then these values of $\alpha$ make inequality \eqref{2.00} valid. To simplify the notations, we set $y=\lambda_{n+1}/\Lambda_{n+1}$ and substitute the upper bound of $\mu_n$ in \eqref{3.14} to see that it suffices to find the values of $0 < \alpha <1$ such that for $0 < y \leq 1$ (note that $y=1$ corresponds to the case $n=1$),
\begin{align*}
  \left ( \frac { \left ( 1 -\left (1-y \right)^{\alpha} \right )^{q/(q-1)}}{y^{1/(q-1)}}+a^{1/(q-1)} \right )^{q-1} \left (\frac {1-y}{a}+ \left(\frac {p-1}{\alpha p} \right )^{q}y \right ) \leq 1.
\end{align*}
   Once again by considering the Taylor expansions, we see that we must take $a=\alpha^q q^{q-1}$ and the above inequality then becomes exactly inequality (3.2) in \cite{G8}. It then follows from the argument in \cite{G8} that this leads to a proof of Theorem \ref{thm2}.

   For the second proof, we set $a_n, b_n$ and $U_p$ as in \eqref{3.13'} and apply \eqref{recurrent} to see that inequality \eqref{2.00} is valid provided that we can define a positive sequence $( \mu_n )_{n \geq 1}$ such that $\mu_1=1, \mu_2=1-U^{-p}_p$ and for $n \geq 2$ (where $\Lambda_0=0$),
\begin{align*}
 \mu_{n+1}=\frac {\left ( 1- \left ( \Lambda_{n-1}/\Lambda_n \right )^{\alpha} \right )^{-p}}{\left (1+\frac {\Lambda_{n-1}}{\lambda_n} \cdot \frac{\lambda_{n-1}}{\Lambda_{n-1}}\left ( 1- \left ( \Lambda_{n-2}/\Lambda_{n-1} \right )^{\alpha} \right )^{-p/(p-1)}\mu^{-1/(p-1)}_n \right )^{p-1}}-\left(\frac {p-1}{\alpha p} \right )^{p}.
\end{align*}
  
  If there exists values of $0< \alpha < 1$ such that for all $n \geq 2$,
\begin{align}
\label{3.17}
  \mu_n \geq \frac {(\frac {\lambda_{n-1}}{\Lambda_{n-1}})^{p-1}}{a^{p-1}\left(1-(1-\frac {\lambda_{n-1}}{\Lambda_{n-1}})^{\alpha} \right )^{p}},
\end{align}
   where $a$ is a positive constant to be determined in what follows, then these values of $\alpha$ make inequality \eqref{2.00} valid. To simplify the notations, we set $y=\lambda_{n}/\Lambda_{n}$ and substitute the lower bound of $\mu_n$ in \eqref{3.17} to see that it suffices to find the values of $0 < \alpha <1$ such that for $0 < y \leq 1$ (note that $y=1$ corresponds to the case $n=2$),
\begin{align*}
  (y+a(1-y))^{1-p} \geq \left ( \frac {p-1}{\alpha p} \right )^{p}y\left ( \frac {1-(1-y)^{\alpha}}{y}\right )^p+a^{1-p}.
\end{align*}
   It is then easy to see that on setting $a=p/(p-1)$, one obtains inequality (3.2) in \cite{G8} again. This gives another proof of Theorem \ref{thm2}.
   
\section{Proofs of Theorems \ref{thm1.4}-\ref{thm1.5} }
\label{sec 5} \setcounter{equation}{0}

   In this section, we use the notations given in \eqref{1.00} and \eqref{1.01} and we note that in order to prove Theorem \ref{Cartimprovement}-\ref{thm1.5}, one may replace the infinite sums by finites sums from $1$ to $N$, with $N \geq 1$ an arbitrary integer. We shall hence assume all the sums in Theorem \ref{Cartimprovement}-\ref{thm1.5} are finite sums from $1$ to $N$ in what follows. Let $n \geq 2$ and let $( a_n )_{n \geq 1}, ( b_n )_{n \geq 1}$ be two positive sequences. We let
\begin{align*}
  \tilde{A}_n=\sum^n_{k=1}\frac {b_kx_k}{a_n}, \quad \tilde{A}^T_{n,N}=b_n\sum^{N}_{k=n}\frac {x_k}{a_k}.
\end{align*}
 We assume that $x_1, \ldots, x_{n-1}$ are fixed real numbers, not necessary positive, such that $\tilde{A}_{n-1} \geq 0$. We regard $\tilde{A}_n$ as a function of $x_n \geq -a_{n-1}\tilde{A}_{n-1}/b_n$ (so that $\tilde{A}_n\geq 0$) only and let
\begin{align*}
g(x_n):=\mu_n \tilde{A}^p_n-\frac {b_n}{a_n} x_n \tilde{A}^{p-1}_n,
\end{align*}
   where $\mu_n<1$ is a constant.
   
   We want to find the maximal value of $g(x_n)$ for $x_n \geq -a_{n-1}\tilde{A}_{n-1}/b_n$. On setting $g'(x_n)=0$, we obtain
\begin{align}
\label{2.1}
   \tilde{A}_n=\frac {b_n}{a_n}\frac {p-1}{p\mu_n-1}x_n.
\end{align}
   We assume $\mu_n \neq 1/p$ for the moment. Using the relation
\begin{align*}
    \tilde{A}_n=\frac {a_{n-1}}{a_n}\tilde{A}_{n-1}+\frac {b_n}{a_n}x_n,
\end{align*}
    we solve $x_n$ to be
\begin{align*}
    x_n=\frac {a_{n-1}}{b_n}\frac {p\mu_n-1}{p(1-\mu_n)}\tilde{A}_{n-1}.
\end{align*}

    Note that the above value of $x_n$ is $> -a_{n-1}\tilde{A}_{n-1}/b_n$. At this value of $x_n$, it is easy to see that
\begin{align}
\label{2.02}
   g''(x_n)=p\left (\frac {b_n}{a_n}\right )^2\tilde{A}^{p-2}_n(\mu_n-1)<0.
\end{align}
   It follows that for $x_n \geq -a_{n-1}\tilde{A}_{n-1}/b_n$, we have
\begin{align*}
    g(x_n) \leq g \left (\frac {a_{n-1}}{b_n}\frac {p\mu_n-1}{p(1-\mu_n)}\tilde{A}_{n-1} \right ).
\end{align*}
   After simplification, the above inequality yields
\begin{align}
\label{2.3}
    \frac {pa_n}{b_n}\mu_n\tilde{A}^p_n- (1-\mu_n)^{1-p}\left (1-\frac 1p \right )^{p-1} \frac {a_n}{b_n}\left (\frac {a_{n-1}}{a_n} \right )^p\tilde{A}^{p}_{n-1} \leq px_n\tilde{A}^{p-1}_n.
\end{align}
   We note that the above inequality continues to hold when $n=1$ with $\mu_1 \leq 1$ and $\tilde{A}_0=0$. It is also easy to check that the above inequality is also valid for $\mu_n=1/p$ when $x_n \geq 0$. Now we let $N \geq 1$ be an integer. Summing the above inequality from $n=1$ to $n=N$ yields
\begin{align*}
   \frac {pa_N}{b_N}\mu_N\tilde{A}^p_{N}+\sum^{N-1}_{n=1} \left (\frac {pa_n}{b_n}\mu_n- (1-\mu_{n+1})^{1-p}\left (1-\frac 1p \right )^{p-1} \left (\frac {a_{n}}{a_{n+1}} \right )^p\frac {a_{n+1}}{b_{n+1}} \right )\tilde{A}^{p}_{n} \leq p\sum^N_{n=1}x_n\tilde{A}^{p-1}_n.
\end{align*}
    We now set $a_n=\Lambda_n, b_n=\lambda_n$ in the above inequality to see that, in order for inequality \eqref{1.40} to hold, it suffices to find a sequence $\{ \mu_n \}$ with $\mu_1=1, \mu_n<1$ when $n \geq 2$ such that
\begin{align*}
   \frac {p\Lambda_n}{\lambda_n}\mu_n- (1-\mu_{n+1})^{1-p}\left (1-\frac 1p \right )^{p-1} \left (\frac {\Lambda_{n}}{\Lambda_{n+1}} \right )^p\frac {\Lambda_{n+1}}{\lambda_{n+1}} \geq p-L.
\end{align*}
   One obtains Theorem \ref{Cartimprovement} on choosing
\begin{align*}
   \mu_n=\frac 1p+(1-\frac 1p)\frac {\lambda_n}{\Lambda_n}.
\end{align*}
   By setting $a_n=\Lambda_n, b_n=\lambda_n$, we can rewrite \eqref{2.3} as
\begin{align*}
    p\Lambda^{1+p-c}_n\mu_nA^p_n- (1-\mu_n)^{1-p}\left (1-\frac 1p \right )^{p-1} \Lambda^{1+p-c}_n\left (\frac {\Lambda_{n-1}}{\Lambda_n} \right )^pA^{p}_{n-1} \leq p\lambda_n\Lambda^{p-c}_na_nA^{p-1}_n.
\end{align*}
   We note that the above inequality continues to hold when $n=1$ with $\mu_1 \leq 1$ and $A_0=0$. Again we set $\mu_1=1$ and let $N \geq 1$ be an integer. Summing the above inequality from $n=1$ to $n=N$ yields
\begin{align}
\label{2.005}
 & \frac{p\Lambda_N}{\lambda_N}\mu_N\lambda_N\Lambda^{p-c}_NA^p_{N}+\sum^{N-1}_{n=1} \left (\frac {p\Lambda_n}{\lambda_n}\mu_n- (1-\mu_{n+1})^{1-p}\left (1-\frac 1p \right )^{p-1} \left (\frac {\Lambda_{n}}{\Lambda_{n+1}} \right )^c\frac {\Lambda_{n+1}}{\lambda_{n}} \right )\lambda_n\Lambda^{p-c}_nA^{p}_{n} \\
& \leq p\sum^N_{n=1}\lambda_n\Lambda^{p-c}_na_nA^{p-1}_n. \nonumber
\end{align}
   We now choose $\mu_n$ to be
\begin{align*}
   \mu_{n}=1-(1-\frac 1p)\left (1-\frac {\lambda_n}{\Lambda_n} \right )^{\frac {c-1}{p-1}}.
\end{align*}
   This choice is made so that 
\begin{align*}
     (1-\mu_{n+1})^{1-p}\left (1-\frac 1p \right )^{p-1} \left (\frac {\Lambda_{n}}{\Lambda_{n+1}} \right )^c\frac {\Lambda_{n+1}}{\lambda_{n}}=\frac {\Lambda_n}{\lambda_n}.
\end{align*}
   It is readily checked that for the so chosen $\mu_n$, when $1< c \leq p$,
\begin{align*}
   &  \frac {p\Lambda_n}{\lambda_n}\mu_n- (1-\mu_{n+1})^{1-p}\left (1-\frac 1p \right )^{p-1} \left (\frac {\Lambda_{n}}{\Lambda_{n+1}} \right )^c\frac {\Lambda_{n+1}}{\lambda_{n}} \\
= & (p-1)\frac {\Lambda_n}{\lambda_n}-(p-1)\frac {\Lambda_n}{\lambda_n}\left (1-\frac {\lambda_n}{\Lambda_n} \right )^{\frac {c-1}{p-1}} \geq c-1.
\end{align*}
    It is easy to see that inequality \eqref{1.8} follows from this.

   We note that inequality \eqref{2.005} continues to hold with $\Lambda_n$ being replaced by $\Lambda^{*}_n$. We now choose $\mu_n=1/p$ so that when $0 \leq c<1$,
\begin{align*}
   &  \frac {p\Lambda^*_n}{\lambda_n}\mu_n- (1-\mu_{n+1})^{1-p}\left (1-\frac 1p \right )^{p-1} \left (\frac {\Lambda^*_{n}}{\Lambda^*_{n+1}} \right )^c\frac {\Lambda^*_{n+1}}{\lambda_{n}} \\
= & \frac {\Lambda^*_n}{\lambda_n}\left (1-\left (1-\frac {\lambda^*_{n}}{\Lambda^*_{n}} \right )^{1-c} \right) \geq 1-c.
\end{align*}
  One readily checks that inequality \eqref{1.10} follows from this. As inequalities \eqref{1.90} and \eqref{1.11} are equivalent to inequalities \eqref{1.10} and \eqref{1.8}, respectively, this completes the proof of Theorem \ref{thm1.5}.

 Now, we let $1 \leq n \leq N-1$ and we assume that $x_{n+1}, \ldots, x_{N}$ are fixed real numbers, not necessary positive, such that $\tilde{A}^T_{n+1,N} \geq 0$. We regard $\tilde{A}^T_{n,N}$ as a function of $x_n \geq -a_{n}\tilde{A}^T_{n+1,N}/b_{n+1}$ (so that $\tilde{A}^T_{n,N}\geq 0$) only and let
\begin{align*}
 h(x_n):=\mu_n({\tilde{A}^T_{n,N}})^p-\frac {b_n}{a_n}x_n ({\tilde{A}^T_{n,N}})^{p-1},
\end{align*}
   where $\mu_n<1$ is a constant.
   
   We want to find the maximal value of $h(x_n)$ for $x_n \geq -a_{n}\tilde{A}^T_{n+1,N}/b_{n+1}$. On setting $h'(x_n)=0$, we obtain a relation given by \eqref{2.1} with $\tilde{A}_n$ being replaced by $\tilde{A}^T_{n,N}$. Again we assume $\mu \neq 1/p$ first and using the relation
\begin{align*}
    \tilde{A}^T_{n,N}=\frac {b_{n}}{b_{n+1}}\tilde{A}^T_{n+1,N}+\frac {b_n}{a_n}x_n,
\end{align*}
    we solve $x_n$ to be
\begin{align*}
    x_n=\frac {a_{n}}{b_{n+1}}\frac {p\mu_n-1}{p(1-\mu_n)}\tilde{A}^T_{n+1,N}.
\end{align*}

    Note that the above value of $x_n$ is $> -a_{n}\tilde{A}^T_{n+1,N}/b_{n+1}$. At this value of $x_n$, it is easy to see that $h''(x_n)$ has an expression similar to that given in \eqref{2.02}, except that one replaces $\tilde{A}_{n}$ by $\tilde{A}^T_{n,N}$ there. It follows that for $x_n \geq -a_{n}\tilde{A}^T_{n+1,N}/b_{n+1}$, we have
\begin{align*}
    h(x_n) \leq h \left (\frac {a_{n}}{b_{n+1}}\frac {p\mu_n-1}{p(1-\mu_n)}\tilde{A}^T_{n+1,N} \right ).
\end{align*}
   After simplification, the above inequality yields
\begin{align*}
    \frac {pa_n}{b_n}\mu_n(\tilde{A}^T_{n,N})^p- (1-\mu_n)^{1-p}\left (1-\frac 1p \right )^{p-1} \frac {a_n}{b_n}\left (\frac {b_{n}}{b_{n+1}} \right )^p(\tilde{A}^{T}_{n+1,N})^p \leq px_n(\tilde{A}^{T}_{n,N})^{p-1}.
\end{align*}
   We note that the above inequality continues to hold when $n=N$ with $\mu_N \leq 1$ and $\tilde{A}^T_{N+1,N}=0$. It is also easy to check that the above inequality is also valid for $\mu_n=1/p$ when $x_n \geq 0$. Summing the above inequality from $n=1$ to $n=N$ (with $\mu_N \leq 1$) yields
\begin{align*}
 &  \frac {pa_1}{b_1}\mu_1(\tilde{A}^T_{1,N})^p+\sum^{N}_{n=2} \left (\frac {pa_n}{b_n}\mu_n- (1-\mu_{n-1})^{1-p}\left (1-\frac 1p \right )^{p-1} \left (\frac {b_{n-1}}{b_{n}} \right )^p\frac {a_{n-1}}{b_{n-1}} \right )(\tilde{A}^{T}_{n,N})^p \\
 \leq & p\sum^N_{n=1}x_n(\tilde{A}^{T}_{n,N})^{p-1}.  \nonumber
\end{align*}
    We now set $a_n=\Lambda_n, b_n=\lambda_n$ in the above inequality to see that, in order for inequality \eqref{1.07} to hold, it suffices to find a sequence $( \mu_n )_{n \geq 1}$ with $\mu_n<1$ such that 
\begin{align}
\label{2.9}
 &  \frac {p\Lambda_n}{\lambda_n}\mu_n- (1-\mu_{n-1})^{1-p}\left (1-\frac 1p \right )^{p-1} \left (\frac {\lambda_{n-1}}{\lambda_{n}} \right )^p\frac {\Lambda_{n-1}}{\lambda_{n-1}} \geq p-(p-1)L, \quad n \geq 2; \\
\label{2.9''}
 & \frac {p\Lambda_1}{\lambda_1}\mu_1 \geq p-(p-1)L. 
\end{align}
   We now set for $n\geq 1$,
\begin{align*}
   \mu_n=1-(1-\frac 1p)\frac {\lambda_n}{\lambda_{n+1}},
\end{align*}
    so that
\begin{align*}
    (1-\mu_{n-1})^{1-p}\left (1-\frac 1p \right )^{p-1} \left (\frac {\lambda_{n-1}}{\lambda_{n}} \right )^p\frac {\Lambda_{n-1}}{\lambda_{n-1}}=\frac {\Lambda_{n-1}}{\lambda_{n}}.
\end{align*}
    It is easy to see that for the so chosen $\mu_n$, inequalities \eqref{2.9} and \eqref{2.9''} are implied by condition \eqref{022}. This completes the proof of Theorem \ref{thm1.4}.

\section{An inequality of Hardy, Littlewood and  P\'{o}lya}
\label{sec 3} \setcounter{equation}{0}
   Our goal in this section is to study inequality \eqref{1} with $C_p=(p/(1-p))^p$ (we remark here that it is shown in \cite{G9} that this value of $C_p$ is best possible) and \eqref{1.9}. In this section we let $1/3 \leq p < 1$. We first note that in order to establish inequality \eqref{1}, it suffices to prove the following inequality:
\begin{align}
\label{2}
  \sum^{N}_{n=1}\Big( \frac 1{n} \sum^{N}_{k=n}x_k \Big
  )^p \geq \left(\frac p{1-p} \right )^p \sum^{N}_{n=1}x^p_n,
\end{align}
   where $N \geq 1$ is an arbitrary integer. 

   By setting 
\begin{align*}
   y_n= \frac 1{n} \sum^{N}_{k=n}x_k, \quad 1\leq n \leq N,
\end{align*}
   we can recast inequality \eqref{2} as (with $y_{N+1}=0$)
\begin{align*}
  \sum^{N}_{n=1}y^p_n \geq \left(\frac p{1-p} \right )^p \sum^{N}_{n=1}(n y_n-(n+1)y_{n+1})^p.
\end{align*}
   Here we require that $y_n \geq 0$ and $y_n \geq (n+1)y_{n+1}/n$. On setting $b_{n}=y_{n}, \lambda_n=\frac 1n, \Lambda_n=1, q=p$ in the function $f(x)$ defined on \cite[p. 843]{G5} and going through the argument there, it is easy to see that for real numbers $\mu_n>n^p$, we have
\begin{align*}
  (n y_n-(n+1)y_{n+1})^p \leq \mu_n y^p_n-(n+1)^p\left (n^{p/(p-1)}\mu^{1/(1-p)}_n-1 \right )^{1-p}y^p_{n+1}.
\end{align*}
   Summing the above inequalities from $1$ to $N$, we obtain
\begin{align*}
  \sum^{N}_{n=1}(n y_n-(n+1)y_{n+1})^p \leq \mu_1 y^p_1+\sum^{N-1}_{n=1}\left(\mu_{n+1}-(n+1)^p\left (n^{p/(p-1)}\mu^{1/(1-p)}_n-1 \right )^{1-p} \right )y^p_{n+1}.
\end{align*}
   We now set
\begin{align*}
  \mu_1=\mu_{n+1}-(n+1)^p\left (n^{p/(p-1)}\mu^{1/(1-p)}_n-1 \right )^{1-p} =\left ( \frac {1-p}{p} \right )^p.
\end{align*}
   Thus, in order for inequality \eqref{2} to hold for some $1/3\leq p <1$, it suffices to show that the sequence $( \mu_n )_{n \geq 1}$ defined above satisfies $\mu_n >n^p$ for this $p$. We now proceed to see for what values of these $p$, we can have $\mu_n \geq an+b$ for all $n \geq n_0$ with $n_0$ some integer $\geq 1$, with the constants $a, b$ to be determined in what follows. Assuming $\mu_n \geq an+b$ for some $n \geq 1$, then using the definition of $\mu_{n+1}$,  we see that it suffices to have
\begin{align*}
  (n+1)^p\left (n^{p/(p-1)}(an+b)^{1/(1-p)}-1 \right )^{1-p}+\left ( \frac {1-p}{p} \right )^p \geq a(n+1)+b.
\end{align*}
   Once again by using Taylor expansions, we see that the optimal values of $a, b$ are $(p/(1-p))^{1-p}, (1/p-1)^p/2$, respectively. With these values of $a,b$ in the above inequality, we see that it becomes inequality (2.1) in \cite{G9} with $y=1/n, t=p/(1-p)$ there and hence is valid by \cite[Lemma 2.1]{G9}. Note that this process also shows that $\mu_{n} > n^p$ for $n \geq n_0$. This now proves Theorem \ref{thm3.2}.

   In the rest of the section, we study inequality \eqref{1.9}, the equivalence of inequality \eqref{1}. Once again, it suffices to consider the following inequality:
\begin{align}
\label{3.5'}
  \sum^{N}_{n=1}\Big(\sum^{n}_{k=1}\frac {x_k}{k} \Big
  )^{q} \leq  \Big ( \frac {p}{1-p} \Big )^q \sum^{N}_{n=1}x^{q}_n,
\end{align}
  where $N \geq 1$ is an arbitrary integer.

  For $1 \leq n \leq N$, we set 
\begin{align*}
  y_n=\sum^{n}_{k=1}\frac {x_k}{k}.
\end{align*}
  We further set $y_0=0$ so that for $1 \leq n \leq N$,
\begin{align*}
  x_n=n(y_n-y_{n-1}).
\end{align*}
  We can thus recast inequality \eqref{3.5'} as
\begin{align*}
  \sum^{N}_{n=1}y^q_n \leq \left ( \frac p{1-p} \right )^q \sum^N_{n=1}(ny_n-ny_{n-1})^q,
\end{align*}
   where we assume $y_n > y_{n-1}>0$ for all $n \geq 2$.

   Note that our argument for the function $f(x)$ defined in \eqref{2.2'} is still valid when $p<0$, provided that we assume $0< x < \alpha y_n/\beta$. On setting $x=y_{n-1}, \alpha=\beta=n, p=q$ there, we see that for all $n \geq 2$, $y_n > y_{n-1}>0, \mu_n \geq 0$,  we have
\begin{align*}
  (ny_n-ny_{n-1})^q+\mu_n y^q_{n-1} \geq \frac {n^q}{(1+n^{q/(q-1)}\mu^{-1/(q-1)}_n)^{q-1}}y^q_n.
\end{align*}
  Together with the observation that $y^q_1 \geq y^q_1$, we see on summing these inequalities from $1$ to $N$ (we may assume $N \geq 2$ here) that we have (with $\mu_{N+1} \geq 0$)
\begin{align*}
 \sum^N_{n=1}(ny_n-ny_{n-1})^q \geq (1-\mu_2)y^q_1+\sum^N_{n=2}\left ( \frac {n^q}{\left (1+n^{q/(q-1)}\mu^{-1/(q-1)}_n \right )^{q-1}}-\mu_{n+1} \right )y^q_n.
\end{align*}
   We now set $\mu_1=0$ and for $n \geq 1$,
\begin{align*}
  \mu_{n+1} =\frac {n^q}{\left (1+n^{q/(q-1)}\mu^{-1/(q-1)}_n \right )^{q-1}}-\left ( \frac 1p-1 \right )^q.
\end{align*}
   It is then easy to see that if the above defined sequence is positive when $n \geq 1$ for some $1/3 \leq p<1$, then inequality \eqref{3.5'} and hence inequality \eqref{1} holds for this value of $p$.

  We now proceed to see for what values of these $p$, we can have $\mu_n \geq an+b$ for all $n \geq n_0$ with $n_0$ some integer $\geq 1$, with the constants $a, b$ to be determined in what follows. Assuming $\mu_{n_0} \geq an_0+b$ for some $n \geq n_0$, then using the definition of $\mu_{n+1}$,  we see that in order for $\mu_n \geq an+b$ to hold for all $n \geq n_0$, it suffices to have, for all $n \geq n_0$,
\begin{align*}
 \left( n^{-p}+(an+b)^{1-p} \right )^{1/(1-p)}-\left (\frac 1p-1 \right )^{p/(p-1)} \geq a(n+1)+b.
\end{align*}
   Once again by using Taylor expansions, we see that the optimal value of $a$ is $(1/p-1)^{1/(p-1)}$. Substituting this value of $a$ in the above inequality, by setting $y=1/n, c=b/a$, we can recast the above inequality as
\begin{align*}
  f_{c,p}(y):=\left (\frac 1p-1 \right )y+(1+cy)^{1-p} -\left (1+ \left (c+\frac 1p \right )y \right )^{1-p} \geq 0.
\end{align*}
  One checks that $f_{c,p}(0)=f'_{c,p}(0)=0$. If we assume $c > -1/(2p)$, then $f''_{c,p}(0) > 0$. It follows that the equation $f''_{c,p}(y)=0$ has at most one root in the interval $(0, 1)$. One then deduces easily that in order for $f_{c,p}(y) \geq 0$ for all $0 \leq y \leq 1/n_0$, it suffices to have $f_{c,p}(1/n_0) \geq 0$. This combined with our discussions above completes the proof of Theorem \ref{thm1.9}.


\end{document}